\documentclass[12pt,twoside]{article}

\usepackage{amssymb}
\usepackage{latexsym}

\newtheorem{example}{Example}[section]

\newtheorem{proposition}[example]{Proposition}

\newtheorem{lemma}[example]{Lemma}
 
\def\Proof{\noindent \it Proof -- \rm}                                                 
\def\qed{\hspace{3.5mm} \hfill \vbox{\hrule height 3pt depth 2 pt width 2mm}
\bigskip}

\def\bo{{\small \,\Box}}
\def\DD{{\sf D}}
\def\EE{{\sf E}}
\def\<{\langle}
\def\>{\rangle}
\def\R{{\mathbb R}}
\def\C{{\mathbb C}}
\def\Z{{\mathbb Z}}
\def\N{{\mathbb N}}
\def\tr{{\rm tr\,}}
\def\gl{{\mathfrak gl}}
\def\glchap{\widehat\gl}
\def\diag{{\rm diag\,}}
\def\Hf{{\rm Hf\,}}
\def\SG{{\mathfrak S}}
\def\H{{\sf H}}
\def\t{{\bf t}}
\def\A{{\sf A}}
\def\aa{{\sf a}}
\def\ch{{\rm ch\,}}
\def\CC{{\cal C}}

\title{Vertex operators and the class algebras of symmetric groups}

\author{Alain {\sc Lascoux} and Jean-Yves {\sc Thibon}}

\date{{\footnotesize Dedicated to the memory of Sergei Kerov}}

\begin{document}

\maketitle

\begin{abstract}
We exhibit a vertex operator which implements multiplication by
power-sums of  Jucys-Murphy elements in the centers of the group algebras
of all symmetric groups simultaneously. The coefficients of this operator
generate a representation of ${\cal W}_{1+\infty}$, to which operators
multiplying by normalized conjugacy classes are also shown to belong.
A new derivation
of such operators based on matrix integrals is proposed, and our vertex
operator is used to give an alternative approach to the polynomial
functions on Young diagrams introduced by Kerov and Olshanski.
\end{abstract}

\section{Introduction}

Convolution of central functions, or
multiplication in the center of the group algebra of the
symmetric group $\SG_n$ can be realized by means of differential
operators acting on symmetric functions. This is due to the
existence of the Frobenius map, which sends a permutation $\sigma$
of cycle type $\alpha$ to the product of power sums $p_\alpha$.
Since the power sums are algebraically independent, any
linear operator on the space $Sym_n$ of homogeneous symmetric
functions of degree $n$ can be realized as a differential operator
in the variables $p_k$, and the above statement is therefore
trivial. More interesting is the existence of infinite order
differential operators  implementing simultaneously
for all symmetric groups
the multiplication by families of elements $\eta_n$ in the
center $Z\SG_n$ of $\C\SG_n$.

The first example of such an operator has been given by
Goulden \cite{Go}, for the case where $\eta_n=\CC_{(2,1^{n-2})}$
is the sum of all transpositions in $\SG_n$, and
similar operators for $\eta_n=\CC_{(\rho,1^{n-r})}$ (where
$\rho$ is a partition of $r$) have been recently described
by Goupil, Poulalhon and Schaeffer \cite{GPS}.

Convolution of central functions is not the only operation
which presents  stability properties allowing to deal simultaneously
with all symmetric groups. This is also the case of  pointwise
multiplication, which, applied to characters, corresponds to
the tensor product of representations. Here, the stability
properties are best explained in terms of vertex operators \cite{T,STW1},
which can also produce stable formulas for outer or inner plethysms,
character values or branching multiplicities \cite{CT,STW2}.

One might therefore expect that vertex operators play a role
in the simultaneous description of the centers of all
symmetric group algebras. A further hint that this should
 be the case is the observation by Frenkel and Wang \cite{FW} that
the commutators of Goulden's operator with the operators
$\alpha_{-k}=$
``multiplication by $p_k$'' and their adjoints $\alpha_k$, generate
a representation of the Virasoro algebra. 

Actually, vertex operators arise when one generalizes Goulden's
formula in another direction. Rather than considering the
sum of all transpositions in $\SG_n$ as the conjugacy class
$\CC_{(2,1^{n-2})}$, one can view it as the sum
$p_1(\Xi_n)=\sum_{i=1}^n\xi_i$ of the Jucys-Murphy elements, and look
for a generating function of the differential operators  
$\DD_k$ implementing
the multiplication by  power sums $p_k(\Xi_n)$ for all $n$.

One finds that the required generating function has a simple
expression in terms of the classical vertex operator
$\Gamma(z_1,z_2)$ describing the Fock space representations
of $\gl_\infty$ and $\glchap_\infty$. Then,
one can see that the brackets  $[\DD_k,\alpha_l]$ generate
a charge 1 representation of the Lie algebra ${\cal W}_{1+\infty}$,
by expressing them in closed form in terms of the standard
generators of this algebra (for $k=1$, this is the result
of Frenkel and Wang).

In Section \ref{sec:appl}, we remark that our vertex operator 
provides  a new approach to the results of Kerov
and Olshanski \cite{KO}. What we  prove is that the coefficients
of the products of power-sums of Jucys-Murphy elements $p_\mu(\Xi_n)$
on the normalized conjugacy classes of \cite{KO} are independent
of $n$, which is equivalent to Proposition 3 of \cite{KO}.

In Section \ref{sec:mat}, we give an alternative derivation of
the  operators of \cite{GPS} in terms of matrix integrals. We
start, following Hanlon, Stanley and Stembridge \cite{HSS}, 
with the observation that the
theory of spherical functions on the cone of positive definite
Hermitian matrices allows one to write generating functions
for connection coefficients as Gaussian integrals over the space
of complex $N\times N$ matrices, $N$ being sufficiently large
(actually, it is the limit $N\rightarrow\infty$ which is relevant,
and we are in fact considering functional integrals).
Then, we combine the generating functions for all $n$, and we
are reduced to the evaluation of similar integrals, but
for a modified Gaussian measure, which can be performed by
means of Wick's formula. On the way, we observe that this
method of calculation can give a direct combinatorial proof
of the generating function of \cite{HSS}, whithout any reference
to spherical functions (this answers a question raised 
in the last section of \cite{HSS}).

Finally, we observe that the results of Kerov and Olshanski
rederived in Section \ref{sec:appl}
imply that the  Goupil-Poulalhon-Schaeffer operators form a linear basis
of the commutative subalgebra of $U({\cal W}_{1+\infty})$
generated by the $\DD_k$.

\section{Notations and background}\label{sec:backg}

\subsection{Symmetric functions}

We denote by $Sym$ the (abstract) algebra of symmetric functions, with
complex coefficients,
and by $Sym(X)=Sym(x_1,\ldots,x_n)$ the algebra of symmetric polynomials
in $n$ variables. The homogeneous component of degree $k$ is denoted
by $Sym_k$. The scalar product on $Sym$ is the standard one, for
which the Schur functions $s_\lambda$ form an orthonormal basis. For
$f\in Sym$, $D_f$ denotes the adjoint of the operator $g\mapsto fg$.
 
For $A=\{a_1,a_2,\ldots\,\}$ we denote by $\sigma_z(A)=\prod_i(1-za_i)^{-1}$
and $\lambda_z(A)=\sigma_{-z}(A)^{-1}$ the generating series of complete
and elementary symmetric functions of $A$. Other notations for
symmetric functions are as in \cite{Mcd}.

\subsection{The Frobenius characteristic map}

The conjugacy class of permutations with cycle type $\mu$ is
denoted by $\CC_\mu$. A central function $f$ on $\SG_n$ is
identified to the element $F=\sum_\sigma f(\sigma)\sigma\in Z\SG_n$.
The Frobenius characteristic map is the linear map $\ch: \C\SG_n\rightarrow
Sym_n$ defined by $\ch(\sigma) =p_\mu$ is $\sigma$ is of cycle type $\mu$.
The structure constants $c_{\alpha\beta}^\gamma$ of $Z\SG_n$
are defined by $\CC_\alpha \CC_\beta=\sum_\gamma c_{\alpha\beta}^\gamma \CC_\gamma$.

The Frobenius map allows one to define a new product $\times$ on each $Sym_n$
by $\ch(F)\times\ch(G)=\ch(FG)$. We extend it to $Sym$ by setting
$u\times v=0$ if $u$ and $v$ are homogeneous of different degrees.
Then, $s_\lambda\times s_\mu = {1\over f^\lambda}\delta_{\lambda\mu}s_\lambda$,
where $f^\lambda$ is the dimension of the representation $\lambda$
of $\SG_n$. We denote by $\Gamma$ the comultiplication dual to $\times$,
that is, $\Gamma(s_\lambda) = {1\over f^\lambda}s_\lambda\otimes s_\lambda$.

\subsection{Jucys-Murphy elements}

The Jucys-Murphy elements of $\SG_n$ are the $n$ sums of
transpositions \cite{Ju,Mu}
\begin{equation}
\xi_j=\xi_{j;n} =\sum_{i<j}(i,j) \,.
\end{equation}
Note that $\xi_1$ is zero, but it is convenient to include
it as well. These elements generate a maximal commutative
algebra $GZ_n$ of $\C\SG_n$
(the Gelfand-Zetlin subalgebra), and the center of $\C\SG_n$
is $Sym(\xi_1,\ldots,\xi_n)$. 
Young's othogonal idempotents $e_{\t,\t}$ ($\t$ a standard
tableau) belong to $GZ_n$, and $\xi_i e_{\t,\t}=c_i(\t)e_{\t,\t}$
where $c_i(\t)$ is the content of the box labelled $i$ in $\t$
(the content of the box in row $k$ and column $l$ of a Young
diagram is defined as $l-k$). The multiset of contents of
a partition $\lambda$ is denoted by $C(\lambda)=\{c_\bo |\bo\in\lambda\}$
(where $\bo$ runs over all boxes in the diagram of $\lambda$).

Jucys has shown that the elementary symmetric function $e_k$ of the
$\xi_i$ is equal to the sum of all permutations having exactly
$n-k$ cycles. One can check that the products 
$$
e_{\bar\alpha}=e_{\alpha_2}e_{\alpha_3} \cdots e_{\alpha_r} \,,
\quad \alpha=(\alpha_1,\alpha_2,\ldots\alpha_r)\vdash n
$$
form a linear basis of $Z\SG_n$. For example, for $n=4$, a basis
is $\{e_0,e_1,e_2,e_3,e_{11}=2\CC_{22}+3\CC_{31}+6\CC_{1111}\}$.
However, in the sequel, we shall rather work with  power sums.

\subsection{The Fock space formalism}

We will also identify $Sym$ with the infinite wedge space ${\cal F}^{(0)}$,
spanned by semi-infinite products $w=v_{i_1}\wedge v_{i_2}\ldots $ 
($i_k\in\Z$) such that $i_1>i_2>\ldots$, and $i_k=1-k$ for $k>>0$.
Such a vector will be denoted by $|\lambda\>$, where
the partition $\lambda$ is defined by $\lambda_k = i_k+k-1$. This space
is the basic representation $L(\Lambda_0)$ of the affine Lie algebra
$\glchap_\infty = {\mathfrak A}_\infty$, the universal central extension
of the Lie algebra of $\Z\times \Z$-matrices with a finite number of
nonzero diagonals. The generators $E_{ij}$ act in a  simple
way of the semi-infinite wedges. For $i\not=j$,
if $v_j$ occurs in a wedge $w$, $E_{ij}$
replaces it with $v_i$, otherwise the result is 0. For $i>0$, $E_{ii}w=w$
if $v_i$ occurs in $w$, and 0 otherwise. For $i\le 0$, the result is
$0$ if $v_i$ occurs in $w$, and $-w$ otherwise.

The Boson-Fermion correspondence
produces differential operators transporting this action on $Sym$
under the linear isomorphism $|\lambda\>\rightarrow s_\lambda$
(see \cite{Kac}, Chap. 14). More generally, the space ${\cal F}^{(m)}$
is spanned by wedges such that $i_k=1-k+m$ for $k>>0$. Their direct sum
(for $m\in\Z$) is called the fermionic Fock space, and
${\cal F}^{(m)}$ is called the charge $m$ sector.

\section{A vertex operator for power sums of Jucys-Murphy elements}\label{sec:vertex}

\subsection{A differential operator for $\SG_n$}

In this section, we will compute the differential operator $\DD^{(n)}$
acting on $Sym_n$ as the generating function
\begin{equation}
F_n(t) = \sum_{k\ge 1}p_k(\Xi_n){t^k\over k!}
=\sum_{i= 1}^n(e^{t\xi_i}-1) \,,
\end{equation}
that is, for $P\in Sym_n$, $\DD^{(n)}P:=\ch(F_n(t))\times P$.

We know that the eigenvalue of $p_k(\Xi_n)$ on the central
idempotent $e_\lambda$ is $p_k(C(\lambda))$. Therefore, the
eigenvalue of $\DD^{(n)}$ on $s_\lambda$ is
\begin{equation}
\sum_{\bo\in\lambda}\left(e^{tc_\bo}-1\right)
=\sum_{\bo\in\lambda}(q^{c_\bo}-1)
\end{equation}
if we set $q=e^t$. This sum is easily evaluated in terms
of the parts of $\lambda$:

\begin{lemma} Let $\lambda$ be a non zero partition of length at most $n$.
Then,
$$ 
\sum_{\bo\in\lambda}q^{c_\bo}=
{q\over q-1}\sum_{i=1}^n\left(q^{\lambda_i-i}-q^{-i}\right)\,.
 $$
\end{lemma}

\Proof The contents of $\lambda$ are the numbers $-i+1,-i+2,\ldots,\lambda_i-1$
for $i=1,\ldots,\ell(\lambda)$.
\qed              

Therefore,
\begin{equation}
\sum_{\bo\in\lambda}(q^{c_\bo}-1)
={q\over q-1}\sum_{i=1}^n\left(q^{\lambda_i-i}-q^{-i}\right)
-\sum_{i=1}^n\lambda_i\,.
\end{equation}

Our first task is to express the operator induced by
$\DD^{(n)}$ on the space of symmetric polynomials
$Sym(x_1,\ldots,x_n)$ in terms of the variables $x_i$.
We set
\begin{equation}
\Delta_n=\prod_{i<j}(x_i-x_j)\quad {\rm and}\quad
\Box_n=(x_1x_2\cdots x_n)^n\,. 
\end{equation}
Let $D_i =x_i{\partial\over\partial x_i}$. We have
\begin{lemma}
$$
{\Box_n\over\Delta_n}\left(
\sum_{i=1}^n q^{D_i} \right) {\Delta_n\over\Box_n}\cdot s_\lambda
=\left(\sum_{i=1}^n q^{\lambda_i-i}\right) s_\lambda \,.
$$
\end{lemma} 

\Proof Multiplication of $s_\lambda$ by $\Delta_n\Box_n^{-1}$ results in  
the determinant $\det(x_j^{\lambda_i-i})$. Applying $\sum_i q^{D_i}$
to this determinant amounts to apply the one-variable operator $q^D$ to
each row of the determinant, and then take the sum. This produces the
same result as applying the operator to each column succesively, 
since both expressions
are equal to the coefficient of $\epsilon$ in
$\det ((1+\epsilon q^D)(x_j^{\lambda_i-i}))$. 
\qed

Therefore, the operator
\begin{equation}
{q\over q-1}{\Box_n\over\Delta_n}\left(
\sum_{i=1}^n q^{D_i}-q^{-i} \right) {\Delta_n\over\Box_n}
-\sum_{i=1}^n D_i
\end{equation}
has the same eigenvalues as $\DD^{(n)}$ on Schur functions $s_\lambda$
in $n$ variables, and must therefore coincide with it. We
can now let $n\rightarrow\infty$, and see that
$\DD^{(n)}$ is  the restriction to $Sym_n$ of the
well-defined limit
\begin{equation}
\DD=\lim_{n\rightarrow\infty}\DD^{(n)}\,.
\end{equation} 
Hence, all symmetric groups can be dealt with simultaneously by
the single operator $\DD$.

\subsection{Bosonization}

The next step is to express $\DD$ in terms of the power sums.
To avoid confusion, we reserve the letter $X$ for finite sets
of variables, and introduce an infinite alphabet $A$ as argument
of our symmetric functions. So, we want to express the
action of $\DD$ on $Sym(A)$ in terms of the operators
\begin{equation} 
\alpha_{-k}=p_k(A)\,,\quad \alpha_k =\alpha_{-k}^\dagger=D_{p_k}
=k{\partial\over\partial p_k(A)}\ \  (k\ge 1) \,.
\end{equation}
This procedure is called {\em bosonization} in the physics
literature
(see, e.g., \cite{AKOS}), for the $\alpha_k$ satisfy the commutation relations
\begin{equation}
[\alpha_j,\alpha_k]=j\delta_{j,-k}
\end{equation}
of  the modes of a  free boson field (a Heisenberg algebra).

To compute the bosonization of $\DD$, we have to calculate
the bi-symmetric kernel
\begin{equation}
K(X,A)= \lambda_{-1}(XA)\DD^{(n)}\sigma_1(XA)
\end{equation}
where $\DD^{(n)}$ acts on functions of $X=\{x_1,\ldots,x_n\}$,
and to express it in the form
\begin{equation}
K(X,A)= \sum_{\mu,\nu}k_{\mu\nu}p_\mu(X)p_\nu(A)\,.
\end{equation}
Then, we will have
\begin{equation} 
\DD=\sum_{\mu,\nu}k_{\mu\nu}p_\mu(A)D_{p_\nu(A)}\,.
\end{equation}
Indeed, writing $\<\,,\,\>$ for the scalar product of $Sym(A)$,
we have $f(X)=\<\sigma_1(XA),f(A)\>$, so that
$$
\DD^{(n)} f(X)=\<\DD^{(n)} \sigma_1(XA),f(A)\>=\<1,D_{\DD^{(n)} \sigma_1(XA)}f(A)\>
=\<\sigma_1(XA),\DD f(A)\>\,.
$$

Let $\nabla_i$ be the partial $q$-derivative with respect to $x_i$,
{\it i.e.}
$$
\nabla_i = {q^{D_i}-1\over (q-1)x_i}\,.
$$
We have, for any function $f(X)=f(x_1,\ldots,x_n)$,
$$
{q\over q-1}\left(\sum_{i=1}^n q^{\lambda_i-i}\right)f(X)
=\left(\sum_{i=1}^n\nabla_i x_i\right)f(X) + \left(\sum_{i=1}^n [i]_{1/q}\right)f(X)
\,.
$$
To apply this to $f(X)=\Delta_n\sigma_1(XA)/\Box_n$, we note that
\begin{equation}
\left(\sum_{i=1}^n\nabla_i x_i\right) {\Delta_n\over\Box_n}
= \left[ {n\over 1-q} + {q(1-q^{-n})\over (1-q)^2}\right] {\Delta_n\over\Box_n} 
\end{equation}
and
\begin{equation}
\left(\sum_{i=1}^n\nabla_i x_i\right) {\Delta_n\over\Box_n}\sigma_1(XA)
=
\left(\sum_{i=1}^n\nabla_i x_i{\Delta_n\over\Box_n}\right)\sigma_1(XA)
+\sum_{i=1}^n {qx_i\Delta^{(i)}\over q^n\Box_n} \nabla_i\sigma_1(XA)
\end{equation}
where 
$$
\Delta^{(i)}=q^{D_i}\Delta_n=\Delta_n A_i(X;q)\,, \quad
 A_i(X;q)=\prod_{j\not =i}{qx_i-x_j\over x_i-x_j}\,.
$$
Hence, setting $\bar\DD^{(n)}=\DD^{(n)}+\EE^{(n)}$, where
$\EE^{(n)}$ is the Euler operator, we have
\begin{eqnarray*}
\bar\DD^{(n)}\sigma_1(XA)&=&
\left[\left({n\over 1-q} + {q(1-q^{-n})\over (1-q)^2}i\right)
+\sum_{i=1}^n[i]_{1/q}
\right]\sigma_1(XA)\\
&&+
q^{1-n}\sum_{i=1}^n x_i A_i(X;q)\nabla_i\sigma_1(XA)
\\
&=&q^{-n}\sum_{i=1}^n A_i(X;q){qx_i\over(q-1)x_i}
\left( {\sigma_{qx_i}(A)\sigma_1(XA)\over\sigma_{x_i}(A)}-\sigma_1(XA)\right)
\\
&=&\sigma_1(XA)
{q^{-n}\over q-1}\sum_{i=1}^n A_i(X;q)(\sigma_{qx_i}(A)\lambda_{-x_i}(A)-1)
\\
&=&\sigma_1(XA)
{q^{-n}\over q-1}\sum_{i=1}^n A_i(X;q)(\sigma_{x_i}((q-1)A)-1)
\\
&=&\sigma_1(XA)
{q^{-n}\over q-1}\sum_{m\ge 1}h_m((q-1)A)\sum_{i=1}^n A_i(X;q)x_i^m
\\
&=&\sigma_1(XA) 
{q^{-n}\over q-1}\sum_{m\ge 1}h_m((q-1)A)q^n{h_m((1-q^{-1})A)\over 1-q^{-1}}
\,.
\end{eqnarray*}
Rewriting this expression  in a more symmetric form,
we obtain the  kernel of $\bar\DD=\DD+\EE$
\begin{equation}
\bar K(X;A)={q\over (q-1)^2}\sum_{m\ge 1}q^{-m}h_m((q-1)A)h_m((q-1)X)\,.
\end{equation}
The bosonization of the Euler operator  being obviously 
$\EE=\sum_{k\ge 1} p_k D_{p_k}$,
we have
\begin{proposition}
The differential operator corresponding to $\sum_{i\ge 1}(q^{\xi_i}-1)$ is
$$
\DD={q\over (q-1)^2}\sum_{m\ge 1}q^{-m}h_m((q-1)A)D_{h_m((q-1)A)}-\EE \,.
$$
\end{proposition}

On this expression, it is clear that $\DD$ can be written
$$\DD= {V_0-1\over (q-1)(1-q^{-1})}-\EE $$
where $V_0$ is the zero mode of the vertex operator
       
\begin{eqnarray*}
V(z;q)&=&\sigma_z((q-1)A)D_{\sigma_{1/z}((1-q^{-1})A)}\\
&=&
\exp\left\{\sum_{k\ge 1}(q^k-1)p_k {z^k\over k}\right\}
\exp\left\{\sum_{l\ge 1}(1-q^{-l})z^{-l}{\partial\over\partial p_l}\right\}\\
&=& :\exp\left\{\sum_{k\not = 0}{(1-q^{-k})z^{-k}\over k}\alpha_k\right\}:\\ 
&=&\sum_{m=-\infty}^\infty V_m(q)z^{-m}\,.
\end{eqnarray*}
This operator satisfies the commutation relations
\begin{equation}
[V,\alpha_k]= z^k(1-q^k)V \quad (k\not = 0)
\end{equation}
so that
\begin{equation}
[V_l,\alpha_k]=(1-q^k)V_{k+l} \quad (k,l\in\Z,\ k\not = 0).
\end{equation}
In particular, 
\begin{equation}\label{eq:Vk} 
V_k = (1-q^k)^{-1}[V_0,\alpha_k] \quad (k\not = 0)
\end{equation}
that is,  all the modes are generated by the action of the
bosonic operators on $V_0$.

\subsection{Class algebras and infinite dimensional Lie algebras}

The vertex operator $V(z;q)$ is well-known to be related to the
Fock space representations of various infinite dimensional
Lie algebras (see {\it e.g.}, \cite{Kac}, Corollary 14.10).
In the notation of \cite{Kac}, $V(z;q)=\Gamma(qz,z)$, and
if  $\hat r_m$ denotes  the representation
of $\glchap_\infty$
in the charge $m$ sector ${\cal F}^{(m)}$ of the fermionic Fock
space ${\cal F}$, one has in particular
\begin{equation}
W(z;q)={1\over 1-q^{-1}}\left(V(z;q)-1\right)
=\sum_{i,j\in\Z}q^i z^j \hat r_0(E_{i,i-j})\,.
\end{equation}
Write $W(z;q)=\sum_k W_k(q) z^{-k}$, so that
\begin{equation}
W_k(q)=\hat r_0 \left(\sum_{i\in \Z} q^i E_{i,i+k}\right)\,.
\end{equation} 
The commutation relations between the operators $W_k(q)$ are
easily determined from the defining relations of $\glchap_\infty$,
which read
\begin{equation}
[E_{ij},E_{kl}] = \delta_{jk}E_{il}-\delta_{li}E_{ki}+\Psi(E_{ij},E_{kl})c\,
\end{equation}
where $c$ is the central charge, and $\Psi$ is the 2-cocycle of
$\gl_\infty$ given by
\begin{equation}
\Psi(E_{ij},E_{ji})=-\Psi(E_{ji},E_{ij})=1\ {\rm if}\ i\le 0,\ j\ge 1
\end{equation} 
and  $\Psi(E_{ij},E_{kl})=0$ in all other cases. One has $\hat r_0(c)=1$,
and a short calculation yields
\begin{equation} 
[W_k(a),W_l(b)]=(b^k-a^l)W_{k+l}(ab)
+\delta_{k,-l}{b^{-l}-a^{-k}\over 1-(ab)^{-1}}\,.
\end{equation}
One recognize that these relations are almost the standard presentation
(in generating function form such as in \cite{KR}, Eq. (2.2.2))
of the Lie algebra usually denoted by ${\hat{\cal D}}$ or ${\cal W}_{1+\infty}$,
the universal central extension of the Lie algebra 
${\cal D}$ of all differential operators on the circle. The generators
of ${\cal D}$ are the $z^kD^n$, where $D=z\partial_z$, and the
corresponding elements of the central extension ${\hat{\cal D}}$
are denoted by $L_k^n$. The cocycle of the central extension
is given by
\begin{equation}
\Phi(z^kf(D),z^lg(D))=\sum_{j\ge 1}kf(-j)g(k-j) \ {\rm if}\ k=-l\ge 0
\end{equation} 
and is 0 in all other cases.
With this at hand, we see that the operators
\begin{equation} 
T_k(q)=-q^{-1}W_k(q)
\end{equation}
satisfy
\begin{equation}
[T_k(a),T_l(b)]
=(a^l-b^k)T_{k+l}(ab)+\delta_{k,-l} {a^{-k}-b^{-l}\over 1-ab}
\end{equation}
which is exactly Eq. (2.2.2) of \cite{KR} with $C=1$. Therefore,
the coefficients $T_{k,n}$ defined by
\begin{equation}
T_k(e^t)=\sum_{n\ge 0} {t^n\over n!}T_{k,n}
\end{equation}
are the images of the $L^k_n$ in a representation $R_0$ of charge 1.

Now, our differential operator $\DD$ reads
\begin{eqnarray*}
\DD &=&{-1 \over 1-q^{-1}} T_0(e^t)-\EE
={1\over t}{-t\over e^{-t}-1}
                       \sum_{l\ge 0}{t^l\over l!}T_{0l}-\EE\\
&=& \sum_{n\ge 1}{t^{n}\over n!}\sum_{k=0}^{n}(-1)^{k-1}
{n \choose k}B_k {T_{0,n+1-k}\over n+1-k}\\
&=& \sum_{n\ge 1}{t^{n}\over n!}\DD_n\,.
\end{eqnarray*}
We have $T_{00}=0$, $T_{01}=-\EE$, and Goulden's operator
is $\DD_1=-{1\over 2}(T_{02}+T_{01})$.
Also, $T_{k,1}=\alpha_k$, $T_{k,2}=2L_k$ where the 
$L_k$ are the charge 1 Virasoro operators considered in \cite{FW}.
Since $[T_0(q),\alpha_k]=(1-q^k)T_k(q)$, we find that
the result of \cite{FW} stating
that the commutators $[\DD_1,\alpha_k]$ generate
a Virasoro algebra can be extended as follows:

\begin{proposition}
The commutators $[\DD_j,\alpha_k]$ generate a charge 1 representation
of ${\cal W}_{1+\infty}$.\qed
\end{proposition}

\subsection{Interpretation of the Virasoro operators}

It would be of interest to have interpretations of the other generators
$T_{ij}$ ($i\not =0$) in terms of natural operations on
$Z\SG=\bigoplus_{n\ge 0}Z\SG_n$. As a step in this direction, we can
propose such an interpretation for the positive part of the Virasoro
algebra.

Let, for $k\ge 1$
\begin{eqnarray}
d'_k&=&\sum_{j\ge 0} p_j D_{p_{j+k}} =\sum_{j\ge 0}\alpha_{-j}\alpha_{j+k}\,,\\
d''_k &=& {1\over 2}
\sum_{1\le i,j\atop i+j=k}\alpha_i\alpha_j\,,\\
d_k &=& d'_k+d''_k\,,
\end{eqnarray}
where $p_0=1$.

Let also $\delta_k$, $\delta'_k$ and $\delta''_k$ be the linear maps
$\C\SG_n\rightarrow\C\SG_{n-k}$ defined on permutations by
\begin{eqnarray}
\delta'_k(\sigma)&=&{1\over (l-1)(l-2)\cdots (l -k)}\sigma^{(k)}\,, \\
\delta''_k(\sigma)&=&{1\over 2 k!}
 \sum_{1\le i,j\atop i+j=k} \,i j\,  \sigma^{(i,j)}
 \,,\\
\delta_k(\sigma)&=&\delta'_k(\sigma)+\delta''_k(\sigma)\,,
\end{eqnarray}
where  $\sigma^{(k)}$ (resp. $\sigma^{(i,j)}$) are defined to be
$0$ if $n,n-1,\ldots,n-k+1$ do not belong to the same cycle
of length $l$ (resp. do not constitute two cycles of lengths
$i$, $j$), and otherwise, $\sigma^{(k)}$ and $\sigma^{(i,j)}$
are the permutations whose cycle decomposition is obtained
by erasing $n,n-1,\ldots,n-k+1$ in the cycle decomposition of $\sigma$.

\begin{proposition}
For $u\in Z\SG_n$, one has
\begin{eqnarray*}
{1\over (n)_k}\ch (\delta_k u) &=& d_k \ch(u)\,,\\
{1\over (n)_k}\ch (\delta'_k u) &=& d'_k \ch(u)\,.
\end{eqnarray*}
\end{proposition}

\Proof A direct calculation.  The numerical factors account for the
orders of conjugacy classes, and could have been suppressed by
using instead the basis
$$
b_\mu={1\over n!}\sum_{\tau\in\SG_n} \tau^{-1}\sigma\tau
$$
where $\sigma$ is any permutation of cycle type $\mu$.  \qed

The operators $d_i$ are the images of the generators $L_i$ of
the Virasoro subalgebra of ${\cal W}_{1+\infty}$ under the
above representation, while the $d'_i$ correspond to the
Witt algebra. Therefore, $L_1$ corresponds to the map considered
in \cite{KOV}, and $L_2$ amounts to erasing $n$ and $n-1$
if they are both in the same cycle, or both fixed points. Similarly,
$L_3$ erases $n,n-1,n-2$ if they are in the same cycle, or constitute
two cycles of lengths $1$ and $2$.

\section{A stability property}\label{sec:appl}

The previous result can be used to express the power-sums of
Jucys-Murphy elements as linear combinations of conjugacy classes.
Indeed, for fixed $n$, the generating function 
$$
J_n(t)=\sum_{k\ge 0} \ch(p_k(\Xi_n)) {t^k\over k!}=\sum_{k\ge 0}J_n^k{t^k\over k!}
$$
is equal to the constant term of 
${1\over (q-1)(1-q^{-1})}(V(z;q)-1)p_1^n$, that  is, to
$$
{1\over (q-1)(1-q^{-1})}\left[
\sum_{k=0}^n{n\choose k} p_1^{n-k}h_k((q-1)A)(1-q^{-1})^k
\right]
$$
which has to be expanded with $q=e^t$ and
$$
h_k((q-1)A)=\sum_{\mu\vdash k} \left[
\prod_i(e^{\mu_i t}-1)\right]
 z_\mu^{-1}p_\mu(A)\,.
$$
This is best accomplished by means of a generating function.
We have
\begin{eqnarray*}
{\cal J}(t) &=& \sum_{n\ge 0}{1\over n!}J_n(t)
=[z^0]{(V(z;q)-1)e^{p_1}\over (q-1)(1-q^{-1})} \\
&=& e^{p_1}\sum_{k\ge 1}{(1-q^{-1})^{k-1}\over k!} {h_k((q-1)A)\over q-1} \\
&=& e^{p_1}\sum_{k\ge 1}{(1-q^{-1})^{k-1}\over k!}
      \sum_{\kappa\vdash k}{p_\kappa(q-1)\over q-1}{p_\kappa(A)\over z_\kappa}\,.
\end{eqnarray*}
For $\kappa\vdash k\ge 1$, let
\begin{equation}
\phi_\kappa(t)={(1-q^{-1})^{k-1}\over k!z_\kappa}{p_\kappa(q-1)\over q-1}|_{q=e^t}
\end{equation}
so that if $\kappa=(1^{k_1}2^{k_2}\cdots\,)$,
\begin{equation}
\phi_\kappa(t)={t^{|\kappa|+\ell(\kappa)-2}\over k! k_1!k_2!\cdots}
\left( 1+O(t)\right)\,,
\end{equation}
and
\begin{eqnarray*}
{\cal J}(t) &=&e^{p_1} \sum_{|\kappa|\ge 1}\phi_\kappa(t)p_\kappa(A)\\
&=&
\sum_{n\ge 0}{1\over n!}\sum_{k=1}^n\sum_{\kappa\vdash k}\phi_\kappa(t)\aa_{\kappa;n}
\end{eqnarray*}
where we have set $\aa_{\kappa;n}=(n)_kp_{\kappa,1^{n-k}}=\ch(a_{\kappa;n})$, where
$a_{\kappa;n}$ are the normalized conjugacy classes defined in \cite{KO}.
Hence, if 
\begin{equation}
\phi_\kappa(t)=\sum_{m\ge |\kappa|+\ell(\kappa)-2}\phi_{\kappa;m}
    {t^m\over m!}
\end{equation}
we obtain
\begin{equation}
J_n^m=\ch( p_m(\Xi_n))
=\sum_{k=1}^{m+1} \sum_{\kappa\vdash k\atop \ell(\kappa)\le m-k+2}
                   \phi_{\kappa;m} \aa_{\kappa;n}\,.
\end{equation}
Observe that the coefficients are independent of $n$. Actually, this
is a special case of a result of Kerov and Olshanski \cite{KO},
which is equivalent to the existence of a similar $n$-independent
expansion of all products of power sums $p_\mu(\Xi_n)$ as linear combinations
of the $a_{\kappa;n}$. This more general result can also be obtained
by the same method, but the expressions of the coefficients
$\phi_{\kappa;\mu}$ become more cumbersome. Instead, we observe
that if we can prove that the coefficients of the expansion
$p_m((\Xi_n)\times \aa_{\kappa;n}$ on the basis $\aa_{\nu;n}$ are
independent of $n$, the
general result will follow by induction.

To prove this, consider the generating function
\begin{equation}
G(t;A,B)=\sum_{n\ge 0}\sum_{m\ge 0} \sum_\kappa 
{1 \over n!m!} t^m 
J_n^m\times \aa_{\kappa;n}(A){p_\kappa(B)\over z_\kappa}\,.
\end{equation}
A calculation similar to the previous one (which is the
case $B=0$) shows that
\begin{equation} 
G(t;A,B)=e^{p_1(A)}\sigma_1(AB)
\sum_{r\ge 1}
{h_r((q-1)A)h_r((1-q^{-1})(B+E))\over (q-1)(1-q^{-1})}
\end{equation} 
where symmetric functions of 
the ``exponential alphabet'' $E$ are defined by $\sigma_t(E)=e^t$
({\it i.e.} $p_1(E)=1$ and $p_k(E)=0$ for $k>1$).
On this expression, it is clear that the coefficient of
${t^m\over m!}{p_\kappa(B)\over z_\kappa}$ in $e^{-p_1(A)}G(t;A,B)$
is a polynomial $\sum_\mu d_{\kappa;m}^\mu p_\mu(A)$, so that
\begin{equation}
p_m(\Xi_n) a_{\kappa;n}=\sum_\mu d_{\kappa;m}^\mu a_{\mu;n}
\end{equation}
the coefficients being independent of $n$, as required.

Here is a table of $n$-independent
expressions of the first power-sums of Jucys-Murphy elements
in terms of normalized conjugacy classes.

\begin{eqnarray*}
p_1(\Xi)&=&{1\over 2}a_2\\
p_2(\Xi)&=& {1\over 3} \,a_3+{1\over 2}\,a_{11}\\
p_3(\Xi)&=& {1\over 4} \,a_4+a_{21}+{1\over 2}\,a_2\\
p_4(\Xi)&=&{1\over 5} \,a_5+{1\over 2} \,a_{22}+a_{31}+{2\over 3}\,a_{111}
+{5\over 3} \,a_3+{1\over 2}\,a_{11}\\
p_5(\Xi)&=&{1\over 6} \,a_6+a_{32}+a_{41}+{5\over 2} \,a_{211}+
{15\over 4}\,a_4+5 \,a_{21}+{1\over 2}\,a_2\\
p_6(\Xi)&=& 
{1\over 7} \,a_7+{1\over 2} \,a_{33}+a_{42}+a_{51}+3\,a_{221}+
3 \,a_{311}+7 \,a_5\\
&&+{5\over 4} \,a_{1111}+{25\over 4} \,a_{22}+15 \,a_{31}+ 
{10\over 3} \,a_{111}+7 \,a_3+{1\over 2}\,a_{11}\\
\end{eqnarray*}                                                                                      

One obtains the expression of each $p_m(\Xi_n)$ from the table
by substituting $[(n-k)!]^{-1}z_{\kappa,1^{n-k}}\CC_{\kappa,1^{n-k}}$
to $a_\kappa$, $\kappa\vdash k$. For example,
\begin{eqnarray*}
p_2(\Xi_n)&=&
{1\over 3}{3\cdot (n-3)!\over (n-3)!}\CC_{3,1^{n-3}}
+{1\over 2} {n!\over (n-2)!}\CC_{1,1,1^{n-2}}\\
&=& \CC_{3,1^{n-3}}+{n\choose 2} \CC_{1^n}\,.
\end{eqnarray*}

\section{A matrix integral approach}\label{sec:mat}

In this section, we will express generating functions
for the coproducts of the elements $\aa_{\rho;n}$ as certain
Gaussian integrals over the space of $N\times N$ complex
matrices. Evaluating these integrals by Wick's formula,
we obtain as a byproduct a new derivation of the differential
operators of \cite{GPS}.

The zonal spherical functions of the Gelfand pair $(GL(N,\C),U(N))$
are known to be expressible in terms of Schur functions (see \cite{Mcd}, Chap. VII,
Sec. 5):
$$
\Omega_\lambda(Z) ={s_\lambda(ZZ^*)\over s_\lambda(N)}\,.
$$
As a consequence, we have  closed form evaluations of the matrix
integrals
\begin{equation}\label{eq:matint}
\int_{M_N(\C)} s_\lambda(AZBZ^*) d\nu(Z)
=2^{|\lambda|}h(\lambda)s_\lambda(A)s_\lambda(B)
\end{equation}
where $h(\lambda)$ is the product of the hook-lengths of $\lambda$,
$A$ and $B$ are arbitrary Hermitian matrices, and $d\nu$
is the  Gaussian probability measure
\begin{equation}
d\nu(Z)=(2\pi)^{-N^2}e^{-{1\over 2}\tr(ZZ^*)}dZ\,,\quad
dZ=\prod_{k,l=1}^N dx_{kl}dy_{kl}\,,\quad z_{kl}=x_{kl}+iy_{kl}\,.
\end{equation}

If $|\lambda|=n$, the right-hand side of (\ref{eq:matint}) is
$$
2^n n! \ {s_\lambda(A)s_\lambda(B)\over f_\lambda}
=2^n n! \ \Gamma(s_\lambda)(A\otimes B)
$$
where $\Gamma$ is the comultiplication dual to the $\times$-product,
induced on $Sym$ by the convolution of central functions, and elements
of $Sym\otimes Sym$ are interpreted as functions of tensor product
of (square) matrices.
Therefore, denoting by $u_\lambda^*$ the adjoint of a basis $u_\lambda$,
\begin{eqnarray*}
\Gamma(p_\lambda) 
&=&\sum_{\alpha,\beta}
\<\Gamma(p_\lambda),p_\alpha^*\otimes p_\beta^*\>p_\alpha\otimes p_\beta\\
&=&{1\over (n!)^2}\sum_{\alpha,\beta}
   \<p_\lambda,C_\alpha\times C_\beta\>p_\alpha\otimes p_\beta\\
&=&{1\over n!}\sum_{\alpha,\beta} c_{\alpha\beta}^\lambda p_\alpha\otimes p_\beta
\end{eqnarray*}
(where we have set $C_\alpha=\ch\CC_\alpha$), so that \cite{HSS}
\begin{equation}\label{eq:HSS}
\int_{M_N(\C)} p_\lambda(AZBZ^*)d\nu(Z)
= 2^n\sum_{\alpha,\beta\vdash n} c_{\alpha\beta}^\lambda 
p_\alpha(A)  p_\beta (B) \,.
\end{equation}

We will now form generating functions for the coproducts.
Let $\rho$ be a partition of $r$. Then,
\begin{eqnarray*}
\int_{M_N(\C)}p_\rho(AZBZ^*){[p_1(AZBZ^*)]^{n-r}\over 2^{n-r}}d\nu(Z)
&=&2^r\sum_{\alpha,\beta\vdash r}c^{\rho 1^{n-r}}_{\alpha\beta}
p_\alpha(A)p_\beta(B)\\
& =& 2^rn!\ \Gamma(p_{\rho 1^{n-r}})\,.
\end{eqnarray*}
Therefore, the coproduct of the element $\aa_\rho=\sum_n \aa_{\rho;n}$
is given by
\begin{eqnarray}
\Gamma\left(\sum_{n\ge r}(n)_rp_{\rho 1^{n-r}}\right)
&=&\int_{M_N(\C)}p_\rho(AZBZ^*)e^{{1\over 2}p_1(AZBZ^*)}d\nu(Z)\\
&=&\int_{M_N(\C)}p_\rho(AZBZ^*)d\mu(Z)
\end{eqnarray}
where
\begin{equation}
d\mu(Z)=d\mu_{A,B}(Z)=(2\pi)^{-N^2}e^{-{1\over 2}\tr(ZZ^*-AZBZ^*)}dZ
\end{equation}
is again a Gaussian measure if we assume that the eigenvalues of $A$ and
$B$ are $<1$. Indeed, one can assume that $A=\diag(a_i)$, $B=\diag(b_i)$,
and in this case, $\tr(ZZ^*-AZBZ^*)=\sum_{i,j}(1-a_ib_j)|z_{ij}|^2$.

The total mass of $d\mu$ is
\begin{equation}
{\cal Z}=\int_{M_N(\C)}d\mu(Z)=\prod_{i,j}{1\over 1-a_ib_j}\,.
\end{equation} 
For a function $f$ on $M_N(\C)$, let
\begin{equation} 
\<f\>={1\over {\cal Z}}\int_{M_N(\C)}f(Z)d\mu(Z)
\end{equation}
denote its expectation value. Since $d\mu$ is Gaussian, we can
make use of Wick's formula, which in this context says the following:
if $f_1,f_2,\ldots,f_m$ are $\R$-linear forms on $M_N(\C)$, we have
\begin{eqnarray}
\<f_1\cdots f_{2k-1}\> &=& 0\, \\
\<f_1\cdots f_{2k}\> &=& \Hf (\<f_i f_j\>)
\end{eqnarray}
where the Hafnian of the matrix $(\<f_i f_j\>)$ is defined by
\begin{equation}
 \Hf (\<f_i f_j\>)_{1\le i,j\le 2k} = \sum \prod_{i=1}^k \<f_{l_i} f_{m_j}\>
\end{equation}
the sum being taken over all pairs of $k$-uples
$L=(l_1<l_2<\cdots < l_k)$ and $M$ such that $l_i<m_i$ and $L\cup M=
[1,2k]$.

In the case at hand, the ``propagators'' are given by
\begin{eqnarray}
\<z_{ij}z_{kl}\> &=&0\,, \\
\<z^*_{ij}z^*_{kl}\> &=&0\,, \\
\<z_{ij}z^*_{kl}\> &=& \delta_{il}\delta_{jk}{2\over 1-a_ib_j}\,.
\end{eqnarray}

Let $M=AZBZ^*$, and let $\sigma$ be the  following permutation
of cycle type $\rho$
\begin{equation}\label{eq:sigma}
\sigma=(12\cdots\rho_1)(\rho_1+1,\rho_1+2,\cdots,\rho_1+\rho_2)\cdots(\cdots r)\,.
\end{equation}
Then, our generating function reads
\begin{eqnarray*}
\<p_\rho(M)\>& =&
\sum_{i_1,\ldots,i_r}
M_{i_1,i_{\sigma(1)}}M_{i_2,i_{\sigma(2)}}\cdots M_{i_r,i_{\sigma(r)}}\\
&=&\sum_{i_1,\ldots,i_r\atop j_1,\ldots,j_r}
a_{i_1}b_{j_1}a_{i_2}b_{j_2}\cdots a_{i_r}b_{j_r}
\<z_{i_1j_1}z_{i_2j_2} \cdots z_{i_rj_r} 
z^*_{i_{\sigma(1)}j_1} \cdots z^*_{i_{\sigma(r)}j_r}\> \\
&=&
\sum_{I,J}a_Ib_J\sum_{\tau\in\SG_r}
\prod_{k=1}^r \delta_{i_k,i_{\sigma\tau(k)}}\delta_{j_k,j_{\tau(k)}}
{2\over 1-a_{i_k}b_{j_k}}\\
&=&
2^r\sum_{I,J}\sum_{\tau\in\SG_r} G_{I,J}
\prod_{k=1}^r \delta_{i_k,i_{\sigma\tau(k)}}\delta_{j_k,j_{\tau(k)}}
\end{eqnarray*}
where we have set $G_{i,j}=a_ib_j(1-a_ib_j)^{-1}$ and
$G_{I,J}=\prod_k G_{i_k, j_k}$
Let us regard the multindices $I,J$ as functions $\{1,2,\ldots,r\}\rightarrow
\N^*$. Then, the above product of Kronecker deltas is zero unless
$I$ is constant on the orbits of $\sigma\tau$, and $J$ is constant
on the orbits of $\tau$. To express the final result, we introduce
the following notation. 
Given a permutation $\tau\in\SG_r$ and  a vector $L=(l_1,l_2,\ldots,l_r)$  of 
positive integers, let
\begin{equation}
p_L^\tau = \prod_k\prod_{\gamma\in\, k- {\rm Cycles}(\tau)}p_{l_{\gamma_1}+l_{\gamma_2}
  +\cdots l_{\gamma_k}}
\end{equation}
product on all $k$-cycles
$\gamma=(\gamma_1,\ldots,\gamma_k)$ of $\tau$.
Then, 

\begin{proposition}
As a symmetric function of the eigenvalues of $A$ and $B$,
$$
\<p_\rho(AZBZ^*)\>=2^r\sum_{l_1,\ldots,l_r\ge 1}\sum_{\tau\in\SG_r}
p_L^{\sigma\tau}(A)p_L^{\tau}(B)\,.
$$
where $\sigma$ is defined in (\ref{eq:sigma}).\qed
\end{proposition}

Note that the above calculation provides an answer to a question raised
at the end of \cite{HSS}, namely, to find a direct combinatorial proof
of (\ref{eq:HSS}). Indeed, to obtain the expectations with respect
to $d\nu$ rather than $d\mu$, one just has to replace the propagators
by $G_{i,j}=a_ib_j/2$, in which case the sum over $L$ disappears
(the only remaining term being for $L=(1,1,\ldots,1)$) and one
finds exactly (\ref{eq:HSS}).

Now, if $\rho$ is a reduced partition (no part equal to 1),
\begin{equation}
\sigma_1(A\otimes B)
\<p_\rho(AZBZ^*)\>=2^r z_\rho\sum_{n\ge r}\Gamma(C_{\rho,1^{n-r}})(A\otimes B)
\end{equation}
and in general, the $\times$-multiplication by any symmetric function
$F$ can be implemented by a differential operator as soon as its
coproduct is known in the form
\begin{equation}
\Gamma(F)=\delta\sigma_1\sum f_{\alpha\beta}p_\alpha\otimes p_\beta\,,
\end{equation}
where $\delta$ is the comultiplication defined by
$\delta(p_\mu)=p_\mu\otimes p_\mu$.
Indeed, for any symmetric functions $G,H$,
\begin{eqnarray*}
\<F \times G,H\> &=& \<F, H\times G\> = \<\Gamma(F),H\otimes G\> \\
&=&\sum_{\alpha\beta}f_{\alpha\beta}
 \<\delta\sigma_1 p_\alpha\otimes p_\beta ,H\otimes G\>\\
&=&\sum_{\alpha\beta}f_{\alpha\beta} 
 \<\delta\sigma_1 , D_{p_\alpha}H\otimes D_{p_\beta} G\> \\
&=& \sum_{\alpha\beta}f_{\alpha\beta} 
 \<\sigma_1 , D_{p_\alpha}H\otimes D_{p_\beta} G\> \\ 
&=&\sum_{\alpha\beta}f_{\alpha\beta}\<D_{p_\beta} G, D_{p_\alpha}H\>\\
&=& \left\<\sum_{\alpha\beta}f_{\alpha\beta}p_\alpha D_{p_\beta} G,H\right\>\,.
\end{eqnarray*}

Hence, we recover the result of \cite{GPS}, proving conjectures
of Katriel \cite{Ka}: 
\begin{proposition}
Let $\rho$ be a reduced partition of $r$. 
For any homogeneous symmetric function $G$ of degree $n\ge r$,
$$
C_{\rho,1^{n-r}}\times G = \H_\rho (G)
$$
where $\H_\rho$ is the differential operator
$$
\H_\rho={1\over z_\rho}
\sum_{l_1,\ldots,l_r\ge 1}\sum_{\tau\in\SG_r}
p_L^{\sigma\tau}D_{p_L^{\tau}}\,.
$$
\end{proposition}
\qed

\section{Final remarks}\label{sec:final}

Recall that $p_\nu\times s_\lambda= (\chi^\lambda_\nu/ f^\lambda)
s_\lambda$, where $f^\lambda$ is the number of standard tableaux of
shape $\lambda$. Kerov and Olshanski have shown that for a
partition $\rho$ of $r$, the function
\begin{equation}\label{eq:KO}
f_\rho(\lambda)= (n)_r{\chi^\lambda_{\rho,1^{n-r}}\over f^\lambda}
\end{equation}
is a polynomial in the ``shifted power-sums''
$$
\tilde p_k(\lambda)=\sum_{i\ge 1}(\lambda_i-i)^k-(-i)^k
$$
which are the eigenvalues on $|\lambda\>=s_\lambda$ of the operators
$$
P_k=\hat r_0\left(\sum_{i\in\Z}i^kE_{ii}\right)
$$
of the Fock space representation of $\glchap_\infty$.
(We have proved an equivalent result in Section \ref{sec:appl}.)

A first consequence of this result is that the elements
$\aa_{\rho;n}=\ch(a_{\rho;n})$  ($a_{\rho;n}\in Z\SG_n$ are the
normalized conjugacy classes of \cite{KO})
have
structure constants independent of $n$: there exist
nonnegative integers $g_{\alpha\beta}^\gamma$ such that
\begin{equation}
\aa_{\alpha;n}\times \aa_{\beta;n} = \sum_\gamma g_{\alpha\beta}^\gamma
\aa_{\gamma;n}
\end{equation}
for all $n$. Therefore, the operators $\A_\rho$
implementing simultaneously the multiplication by all
$\aa_{\rho;n}$ also satisfy
\begin{equation}
\A_{\alpha}\A_{\beta} = \sum_\gamma g_{\alpha\beta}^\gamma
\A_{\gamma}
\end{equation}    
so that they form a linear basis of  commutative subalgebra
of $U(\glchap_\infty)$. When $\rho$ is reduced, $\A_{\rho}=z_\rho\H_\rho$.

A second consequence is that these operators actually belong
to the image of $U({\cal W}_{1+\infty})$ under the
representation $R_0$ of Section \ref{sec:vertex}. 
Indeed, (\ref{eq:KO}) show
that $\A_\rho$ is a polynomial in the commuting operators
$P_k$, which are related to the $\DD_k$ of the previous section by
\begin{equation}
P_n=\sum_{k=0}^{n-1}{n\choose k}\DD_k
\end{equation}
where we have set $\DD_0=\EE$.

Kerov and Olshanski also identified the algebra of the $a_{\rho;n}$
to an algebra of differential operators $A_{\rho;N}$ living
in the center of $U(\gl_N)$ for $N\ge n$. Since these operators
commute with the adjoint representation, they preserve the space
of functions on $GL_N(\C)$ which are symmetric
functions of the eigenvalues. The previous considerations show
then that if we set
$$
\tilde p_k(D)=\sum_{i=1}^N D_i^k-(-i)^k
$$
and $f_\rho(\lambda)=\sum_\nu k_{\rho\nu}\tilde p_\nu(\lambda)$,
the restriction of  $A_{\rho;N}$ to invariant functions is given by
\begin{equation}
A_{\rho;N} = {\Box_N\over \Delta_N}
\sum_\nu k_{\rho\nu}\tilde p_\nu(D) {\Delta_N\over \Box_N}\,.
\end{equation}

\bigskip
{\footnotesize\noindent
\sc Institut Gaspard Monge, \\
Universit\'e de Marne-la-Vall\'ee,\\
77454 Marne-la-Vall\'ee cedex,\\
 FRANCE}


\begin{thebibliography}{abc}
%
\bibitem{AKOS}{H. Awata, H. Kubo, S. Odake, and J. Shiraishi}, 
{\it Quantum ${W}_N$ algebras and Macdonald polynomials}, 
Comm. Math. Phys. {\bf 179} (1996), 401--416. 
%
\bibitem{CT}{C. Carr\'e and J.-Y. Thibon},
{Plethysm and vertex operators},
Adv. Appl. Math. {\bf 13} (1992), 390-403. 
%
\bibitem{FW}{I.B. Frenkel and W. Wang}, {Virasoro algebra and wreath
product convolution}, preprint {\tt QA/0006087}.
%
\bibitem{Go}{I. Goulden}, {\it A differential operator for symmetric
functions and the combinatorics of multiplying transpositions},
Trans. Amer. Math. Soc. {\bf 344} (1994), 421--440.
%
\bibitem{GPS}{A. Goupil, D. Poulalhon and G. Schaeffer},
{\it Central characters and conjugacy classes of the symmetric
group}, Proceedings of FPSAC'00 (D. Krob, A.A. Mikhalev and A.V. Mikhalev
eds.), Moscow, June 2000 (Springer), 238--249.
%
\bibitem{HSS}{P.J. Hanlon, R.P. Stanley and J.R. Stembridge}, {\it
Some combinatorial aspects of spectra of normally distributed random
matrices}, Contemp. Math. {\bf 138} (1992), 151--174.
%
\bibitem{Ju}{A. Jucys}, {\it Symmetric polynomials and the center of the
symmetric group rings}, Rep. Math. Phys. {\bf 5} (1974), 107--112.
%
\bibitem{Kac}{V.G. Kac}, {\it Infinite dimensional Lie algebras},
3rd edition, Cambridge, 1990.
\bibitem{Ka}{J. Katriel},{\it The class algebra of the symmetric group}, 
Proc. 10th FPSAC Conference, N. Bergeron and F. Sottile Eds.,
Toronto (1998),401-410.
%
\bibitem{KR}{V.G. Kac and A. Radul}, {\it Quasifinite highest weight modules
over the Lie algebra of differential operators on the circle}, 
Comm. Math. Phys. {\bf 157} (1993), 429--457.
%
\bibitem{KO}{S. Kerov and G. Olshanski}, {\it Polynomial functions on the
set of Young diagrams}, C.R. Acad. Sci. Paris S\'er. I {\bf 319}
(1994), 121--126.
%
\bibitem{KOV}{S. Kerov, G. Olshansi and A. Vershik}, {\it Harmonic
analysis on the infinite symmetric group. A deformation of the
regular representation}, C.R. Acad. Sci. Paris S\'er. I {\bf 316}
(1993), 773--778.
%
\bibitem{Mcd}{I.G. Macdonald}, {\it Symmetric functions and Hall polynomials},
2nd edition, Oxford, 1995.
%
\bibitem{Mu}{G. Murphy}, {\it A new construction of Young's seminormal
representation of the symmetric group}, J. Algebra {\bf 69} (1981), 287--291.
%
\bibitem{OO}{A. Okounkov and G. Olshanski}, {\it Shifted Schur
functions}, St. Petersburg Math. J. {\bf 9} (1998), 239--300.

\bibitem{STW1}{T. Scharf, J.-Y. Thibon and B.G. Wybourne}
 {\it Reduced notation, inner plethysm and the symmetric group},
J.  Phys. A {\bf 24} (1993), 7461--7478.
%
\bibitem{STW2}{T. Scharf, J.-Y. Thibon and B.G. Wybourne} 
{\it Generating functions for stable branching coefficients
of $U(n)\downarrow S_n$, $O(n)\downarrow S_n$ and $O(n)\downarrow S_{n-1}$},
J. Phys. A {\bf 30} (1997), 6963--6975.
%
\bibitem{T}{J.-Y. Thibon},
 {\it Hopf algebras of symmetric functions and tensor products of
symmetric group representations}, 
Internat. J. Alg.  Comp. {\bf 2} (1991), 207--221.                                
%
\end{thebibliography}
\end{document}